\documentclass{amsart}
\usepackage{graphics}
\usepackage{color}
\usepackage{amssymb}
\usepackage{amsmath}
\usepackage[all]{xy}
\theoremstyle{plain}
\newtheorem{theor}{Theorem}[section]
\newtheorem{prop}[theor]{Proposition}
\newtheorem{lemma}[theor]{Lemma}
\newtheorem{cor}[theor]{Corollary}
\theoremstyle{definition}
\newtheorem{de}[theor]{Definition}
\newtheorem{example}[theor]{Example}
\theoremstyle{remark}
\newtheorem{remark}[theor]{Remark}
\def\Chi{{\mathbb X}}

\def\div{{\rm div}}

\def\reg{{\rm reg}}

\def\KK{{\mathbb K}}
\def\ZZ{{\mathbb Z}}

\def\NN{{\mathbb N}}
\def\QQ{{\mathbb Q}}
\def\PP{{\mathbb P}}

\def\Of{{\mathcal{O}}}

\def\Ff{{\mathcal{F}}}

\def\WDiv{\operatorname{WDiv}}

\def\Aut{\operatorname{Aut}}

\def\Cl{\operatorname{Cl}}

\def\GL{{\rm GL}}
\def\Sp{{\rm Sp}}
\def\Pf{{\rm Pf}}

\def\SO{{\rm SO}}

\def\codim{{\rm codim}}
\def\As{{\rm Ass}}
\def\Spec{{\rm Spec\,}}

\def\diag{{\rm diag}}

\def\SL{{\rm SL}}
\def\Ker{{\rm Ker}\,}
\def\Aut{{\rm Aut}}
\def\Imm{{\rm Im}\,}

\newcounter{itemnumber}

\begin{document}
\title [Cox rings and affine varieties]{Cox rings, semigroups and automorphisms of affine algebraic
varieties}
\author [I.V. Arzhantsev and S.A. Gaifullin]{Ivan V. Arzhantsev and Sergey A. Gaifullin}
\address{Department of Higher Algebra, Faculty of Mechanics and Mathematics, Moscow State University,
Leninskie Gory, Moscow, 119991 Russia}
\email{arjantse@mccme.ru}
\email{sgayf@yandex.ru}
\thanks{The first author gratefully acknowledges the support from P.~Deligne 2004 Balzan prize in mathematics}
\begin{abstract}
We study the Cox realization of an affine variety, i.e.,
a canonical representation of a normal affine variety with
finitely generated divisor class group as a quotient of a
factorially graded affine variety by an action of the Neron-Severi
quasitorus. The realization is described explicitly for the quotient space of
a linear action of a finite group. A universal property of this
realization is proved, and some results on the divisor theory of an
abstract semigroup emerging in this context are given. We show
that each automorphism of an affine variety can be lifted to an
automorphism of the Cox ring normalizing the grading. It follows
that the automorphism group of a non-degenerate affine toric variety of
dimension $\geq2$ has infinite dimension. We obtain a wild automorphism of the
three-dimensional quadratic cone that rises to Anick's
automorphism of the polynomial algebra in four variables.
\end{abstract}
\maketitle
\section{Introduction}\label{sec1}
Let $X$ be a normal algebraic variety without nonconstant
invertible functions over an algebraically closed field
$\mathbb{K}$ of characteristic zero. Suppose that the divisor class
group $\mathrm{Cl}(X)$ is free and finitely generated. Denote by
$\mathrm{WDiv}(X)$ the group of Weil divisors on $X$. Let us fix a
subgroup $K\subset \mathrm{WDiv}(X)$ that projects onto
$\mathrm{Cl}(X)$ isomorphically. Following famous D.~Cox's
construction from toric geometry \cite{Cox} (see also \cite[Def.~2.6]{HK}), 
we define the Cox ring (or the total coordinate ring) of the
variety $X$ as
$$
R(X)=\bigoplus_{D\in K}\mathcal{O}(X,D), \ \ \mathrm{where}\ \
\mathcal{O}(X,D)=\{f\in \mathbb{K}(X)\mid \mathrm{div}(f)+D\geq 0
\}.
$$
Multiplication on graded components of $R(X)$ is given by
multiplication of rational functions,
and extends to other elements by distributivity. It is not
difficult to check that the ring $R(X)$ depends on the choice of
the lattice $K$ only up to isomorphism. The ring R(X) is proved to
be factorial, see \cite{BHH} and \cite{EKW}. There is a proof of
this fact based on the concept of graded factoriality \cite{Ar}, that is
uniqueness of the prime decomposition in the
semigroup $\mathrm{Ass}(R(X))$ of association classes of nonzero
homogeneous elements of $R(X)$. This semigroup is
naturally isomorphic to the semigroup $\mathrm{WDiv}(X)_+$ of
effective Weil divisors on $X$, hence it is free.

If the group $\mathrm{Cl}(X)$ is finitely generated but not free,
the definition of the Cox ring $R(X)$ requires more efforts, see
Section \ref{coxring}. Here $R(X)$ is also
factorially graded, but may be not factorial, see
\cite[Example~4.2]{Ar} and Example~\ref{exq8} below.

A quasitorus $N$ is called the {\it Neron-Severi quasitorus} of the variety $X$, if 
an isomorphism between the character group $\mathfrak{X}(N)$ and the group 
$\mathrm{Cl}(X)$ is fixed. Suppose that $R(X)$ is finitely generated.
The quasitorus $N$ acts naturally on the spectrum $\overline{X}$ of
$R(X)$, and $X$ is the categorical quotient for this action. The
quotient morphism $q\colon \overline{X}\rightarrow X$ is said to be
the {\it Cox realization} of the variety $X$. This construction
has already been used in affine algebraic geometry.
For example, the calculation of the Cox ring
allowed to characterize affine $\mathrm{SL}(2)$-embeddings that
admit the structure of a toric variety~\cite{Ga}. In \cite{BaH}, it is shown
that Cox realization leads to a remarkable unified description of
all affine $\mathrm{SL}(2)$-embeddings, which was not known
before. The Cox realization can be defined for a wide class
of normal (not necessarily affine) varieties. If the Cox ring is
finitely generated, the realization results in a combinatorial
description of varieties via the theory of bunched rings.
This approach allows to describe many of geometric properties, see
\cite{BHC} and \cite{Hau}.

We construct explicitly the Cox realization of the quotient space of
a finite linear group. As a corollary, quotient
spaces that are toric varieties are described. In
Section~\ref{sec2}, we prove the following universal property: the
quotient morphism of a factorially graded affine variety by the 
quasitorus action that does not contract divisors can be passed
through the Cox realization of the quotient space.

For an affine variety $X$, the embedding of the algebra of regular
functions $\mathbb{K}[X]$ into $R(X)$ as zero homogeneous
component defines an embedding of the semigroup $\Gamma$ of principle
effective divisors on $X$ into $\mathrm{WDiv}(X)_+$. This 
is the divisor theory for the semigroup $\Gamma$ in the sense of 
\cite[Ch.~III]{BS} and \cite{Sk}. Some problems on Cox rings can be
reformulated and solved in terms of semigroup theory. Conversely,
geometric intuition advises assertions about semigroups that are
of independent interest. One of the aims of this paper is to show
effectiveness of this correspondence. For example, the 
universal property of the Cox realization mentioned above leads to the following
result: an embedding of a semigroup $\Gamma$ into a free semigroup
$S$ under certain conditions can be extended to an embedding of the
divisor theory of $\Gamma$ into $S$ (Theorem \ref{dtheo}). All
definitions and statements concerning semigroups are gathered in
the last section.

One of the main results of the paper \cite{Cox} is a description of
the automorphism group of a complete simplicial toric variety using
homogeneous automorphisms of its Cox ring. Later D.~Buhler
(1996) generalized this result to arbitrary complete toric
varieties. We prove that each automorphism of a normal affine
variety lifts to an automorphism of the Cox ring normalizing the
grading. It is known that the Cox ring of a toric variety is a
polynomial ring \cite{Cox}. This allows to prove that the automorphism
group of a non-degenerate affine toric variety of dimension $\geq 2$ has infinite
dimension.

We say that an automorphism of an affine toric variety $X$ is tame
if it can be represented as a composition of automorphisms lifting 
to elementary automorphisms of the polynomial algebra $R(X)$. 
Other automorphisms are called wild. In Section
\ref{sec4} we give an example of a wild automorphism of 
three-dimensional quadratic cone. Here we do not
use recent results of I.P.~Shestakov and U.~Umirbaev \cite{SU}.
Our automorphism rises to a famous Anick's automorphism of the polynomial algebra, 
see \cite[p.~49, p.~96]{Es} and \cite{Sm}. Anick's automorphism is one
of the candidates to be a wild automorphism of the algebra of polynomials in four
variables. We show that this automorphism can not be decomposed into
a composition of elementary automorphisms preserving the grading.

The authors are grateful to E.B.~Vinberg, who attracted their
attention to the relation between graded factoriality and the divisor
theory of a semigroup, to A.~van den Essen, who referred them to
papers on Anick's automorphism, and to J.~Hausen for fruitful
discussions.

\section{The Cox ring of an affine variety}\label{coxring}
Consider an irreducible normal variety $X$ with finitely generated
divisor class group $\mathrm{Cl(X)}$ and without nonconstant invertible
regular functions. Let $K\subset\mathrm{WDiv}(X)$ be a finitely generated subgroup projecting
surjectively onto $\mathrm{Cl}(X)$. Consider a ring
$$
T_K(X)=\bigoplus_{D\in K}\mathcal{O}(X,D).
$$
Let $K^0\subset K$ be the kernel of the projection
$K\rightarrow\mathrm{Cl}(X)$. Fix some bases
$D_1,\ldots,D_s$ in $K$ and $D_1^0=d_1D_1,\ldots,D_r^0=d_rD_r$ in
$K^0$, $r\leq s$. Let us call a set of rational functions
$\mathcal{F}=\{F_D\in\mathbb{K}(X)^\times:D\in K^0\}$ {\it
coherent} if the principal divisor $\mathrm{div(F_D)}$
coincides with $D$ and $F_{D+D'}=F_DF_{D'}$. It is clear that the
set $\{F_D\}$ is determined by functions $F_{D_i^0},
i=1,\ldots,r$: if $D=a_1D_1+\ldots+a_rD_r^0$, then
$F_D=F_{D_1^0}^{a_1}\ldots F_{D_r^0}^{a_r}$. Let us fix a
coherent set $\mathcal{F}$.

Suppose that $D_1,D_2\in K$ and $D_1-D_2\in K^0$. The map $f\mapsto
F_{D_1-D_2}f$ is an isomorphism of vector spaces
$\mathcal{O}(X,D_1)$ and $\mathcal{O}(X,D_2)$. It is easy to see
that the vector space spanned by vectors $f-F_{D_1-D_2}f$ for all
$D_1,D_2\in K$, $D_1-D_2\in K^0$, $f\in\mathcal{O}(X,D_1)$, is an
ideal $I(K,\mathcal{F})\lhd T_K(X)$. Define the Cox ring of the
variety $X$ as
$$
R_{K,\mathcal{F}}(X)=T_K(X)/I(K,\mathcal{F}).
$$
The ring $R_{K,\mathcal{F}}(X)$ does not depend on the choice of
$K$ and $\mathcal{F}$ up to isomorphism (see, for example,
\cite[Prop.~3.2]{Ar}). Further we denote it by $R(X)$. Since
$K/K^0\cong\mathrm{Cl}(X)$, the ring $R(X)$ has a natural
$\mathrm{Cl}(X)$-grading.
\begin{de}
Suppose $M$ is a finitely generated abelian group and
$R=\bigoplus_{u\in M}R_u$ is an $M$-graded algebra.
\begin{itemize}
\item
A nonzero homogeneous element $a\in R\setminus R^\times$ is called
{\it $M$-irreducible}, if the condition $a=bc$, where $b$ and $c$
are homogeneous, implies that either $b$ or $c$ is invertible.
\item
The graded algebra $R$ is called {\it factorially graded}, if $R$ does not
contain nonconstant homogeneous invertible elements, each
nonconstant homogeneous element can be represented
as a product of $M$-irreducible elements and this representation is
unique up to association and order of factors.
\end{itemize}
\end{de}

Recall that two homogeneous elements $a,b\in R$ are {\it
associated}, if $a=bc$, where $c$ is invertible. The set
$\mathrm{Ass}(R)$ of association classes of non-zero homogeneous
elements of $R$ is a semigroup with multiplication
induced by multiplication in $R$. It is clear that graded factoriality
of an algebra $R$ means that the semigroup $\mathrm{Ass}(R)$
is free. If $R=R(X)$ and $M=\mathrm{Cl}(X)$, then the semigroup
$\mathrm{Ass}(R)$ may be identified with $\mathrm{WDiv}(X)_+$.
Here R(X) contains neither homogeneous zerodivisors nor
nonconstant homogeneous invertible functions, and is
factorially graded.

Recall that a {\it quasitorus} is an affine algebraic group
isomorphic to the direct pro\-duct of a torus and a finite abelian
group. The group of characters of a quasitorus is a finitely
generated abelian group whose rank equals the dimension of
the quasitorus. It is easy to check that each finitely generated
abelian group can be realized as the group of characters of a
quasitorus, and a quasitorus is 
determined by its group of characters up to isomorphism. Define the {\it Neron-Severi
quasitorus} of a variety $X$ as a quasitorus $N$ whose character
group is identified with $\mathrm{Cl}(X)$. The quasitorus
$N$ acts by automorphisms of the algebra $R(X)$: a homogeneous
component $R(X)_u$ is the weight subspace corresponding to the weight
$u$. The subalgebra $R(X)^N$ of invariants of this action coincides
with $R(X)_0=\mathbb{K}[X]$.

Further we assume that the variety $X$ is affine and
$R(X)$ is finitely generated. Consider the spectrum
$\overline{X}=\mathrm{Spec}(R(X))$ and the regular
$N$-action on the affine variety $\overline{X}$. 
Since $X\cong \mathrm{Spec}(R(X)^N)$, the variety $X$
is realized as the $N$-quotient of $\overline{X}$. 
The quotient morphism $q\colon \overline{X}\rightarrow X$ is said to be 
the {\it Cox realization} of the variety $X$.

\begin{prop}\label{nor}
The variety $\overline{X}$ is irreducible and normal.
\end{prop}
\begin{proof}
Let $X^{\mathrm{reg}}$ be the set of smooth points of $X$. It
is an open subset with the complement of codimension
$\geq2$. Fix a finitely generated subgroup $K\subset
\mathrm{WDiv}(X)$ that projects onto $\mathrm{Cl}(X)$ surjectively, and
consider the ring $T_K(X)$. This ring has no homogeneous zerodivisors 
and hence has no zerodivisors at all \cite[Lemma~2.1]{Ar}. 
Every point $x\in X^{\mathrm{reg}}$ has a neighborhood $V$,
where any divisor is principle. The preimage of this
neighborhood in the spectrum of the ring $T_K(X)$ is isomorphic to
$V\times T$, where $T$ is a torus whose character lattice 
is identified with $K$. Therefore the projection of the preimage
$W$ of the variety $X^{\mathrm{reg}}$ is a locally trivial
bundle $p\colon W\rightarrow X^{\mathrm{reg}}$ with fiber $T$.
Geometrically the quotient of $T_K(X)$ over the ideal $I(K,\Ff)$ 
corresponds to the subvariety
$U=q^{-1}(X^{\reg})\subseteq W$. This subvariety intersects each fiber
(isomorphic to $T$) in a subgroup isomorphic to $N$.

It suffices to show that the variety $U$ is irreducible. Indeed,
for dimension reasons $U$ is contained in an irreducible component of
$\overline{X}$ that is the closure of $U$.
Consequently, this component is $N$-invariant. The union of other
components is also $N$-invariant, and if it is non-empty, then there are
homogeneous zerodivisors in $R(X)$.

The subvariety $q^{-1}(V)$ is defined in $p^{-1}(V)=V\times T$ by
equations $t_i^{d_i}=F_{D_i^0}$, where functions $F_{D_i^0}$ are
irreducible in $\KK[X]$. Hence $q^{-1}(V)$ is irreducible and is
contained in one component of $U$. Such components are
$N$-invariant and
cover $U$. Since $X^{\reg}$ is irreducible and the image of
a closed $N$-invariant subset under the morphism $q$ is closed,
the variety $U$ is irreducible.

The proof of the fact that $R(X)$ is integrally closed may be found
in \cite[Prop.~2.7, Rem.~2.8]{Ar}.
\end{proof}

It is easy to check that the algebra $R(X)$ is factorially graded if and only if
each $N$-invariant Weil divisor on $\overline{X}$ is
principle. In particular, $N$-linearized line bundles on
$\overline{X}$ are $N$-linearizations of the trivial one.

The following proposition contains a necessary and sufficient
condition for the quotient morphism of an affine variety by an
action of a quasitorus to be the Cox realization of the quotient space. This
condition emerged in several former papers, see for example
\cite[Rem.~4.2]{AH} and \cite[Th.~2.2]{BaH}.

\begin{prop}\label{prrel}
Let $q\colon\overline{X}\rightarrow X$ be the Cox
realization of an affine variety $X$. Suppose $Q$ is a quasitorus
acting on an irreducible affine variety $Z$, which is factorially graded
with regard to this action. Let $\pi\colon Z\rightarrow Z/\!\!/Q$
be the quotient morphism. Suppose that there exists an isomorphism
$\phi\colon X\rightarrow Z/\!\!/Q$. Then the following conditions
are equivalent:
\begin{enumerate}
\item there exists an open $Q$-invariant subset $U\subseteq Z$ such that
 $\codim_Z(Z\setminus U)\ge 2$, the $Q$-action on $U$ is free
and each fiber of $\pi$ having nonempty intersection with $U$
is a $Q$-orbit;
\item there exist an isomorphism $\mu: Q\to N$ and an isomorphism
$\xi:Z\to\overline{X}$ such that $\xi(gz)=\mu(g)\xi(z)$ for all
$g\in Q$, $z\in Z$, and the following diagram is commutative:
$$
\xymatrix{ & Z \ar[d]^{\pi} \ar[rr]^{\xi} & & \overline{X}
\ar[d]_{q} &
\\
& Z/\!/Q \ar[rr]^{\phi} & & X. & }
$$
\end{enumerate}
\end{prop}

\begin{proof}
For the Cox realization $q:\overline{X}\to X$ one can take
$q^{-1}(X^{\text{reg}})$ as $U$, see the proof of 
Proposition~\ref{nor}.
Conversely, consider the restriction of $\pi$ to $U$:
$$
\pi:U\to U_0\subset Z/\!/Q.
$$
Reducing $U$, we may assume that $U_0$ consists of smooth points
of the variety $Z/\!\!/Q$. For each divisor $D$ on $Z/\!/Q$ its intersection
with $U_0$ is a locally principle divisor, thus it determines a line
bundle $L_D$. The lift of this bundle to $U$ gives a
$Q$-linearization of the trivial bundle on $U$, and hence on $Z$.
Conversely, every $Q$-linearization of the trivial line bundle on
$Z$ comes from a divisor on $U_0$. This defines an isomorphism
between the group of characters of the quasitorus $Q$ and the Picard
group of $U_0$ coinciding with $\Cl(X)$
\cite[Prop.~4.2]{KKV}. The space of invariant sections of the
$Q$-linearized trivial line bundle on $Z$ with respect to a
character $\xi\in\Chi(Q)$ is identified with the subspace of
semiinvariants of weight $-\xi$ in $\KK[Z]$. The
algebra $\KK[Z]$ is a direct sum of subspaces of
$Q$-semiinvariants, so we get a homogeneous isomorphism between
$\KK[Z]$ and $R(X)$.
\end{proof}

Finally we give an example of an affine variety with 
a non-finitely generated Cox ring. If
$Y\subseteq\PP^n$ is an irreducible projectively normal projective
subvariety and $X\subseteq\KK^{n+1}$ is the affine cone over $Y$,
then $R(Y)$ and $R(X)$ are isomorphic as (non-graded) algebras.  
Indeed, for every divisor $D$ on $Y$ the quotient morphism $\pi: X\setminus\{0\}
\xrightarrow{/\!/\KK^{\times}} Y$ determines an isomorphism
between the space $\Of(Y,D)$ and a weight component with respect
to the action of the torus $\KK^\times$ in the space $\Of(X, D')$,
where $D'$ is the image of the class of $D$ in the exact sequence
$0\to\ZZ\to\Cl(Y)\to\Cl(X)\to 0$. But $\Of(X,D')$ is a direct sum
of weight components corresponding to divisors $D$ mapping to
$D'$~\cite[Prop.~4.2]{KKV}.

Let us take as $Y$ the blow up of the projective plane
$\mathbb{P}^2$ at nine generic points. 
By~\cite{Na} (see also \cite[Rem.~3.2]{EKW}),
the ring $R(Y)$ is not finitely generated.


\section{The Cox realization of a quotient space} \label{sec1-1}

Let $V$ be a finite-dimensional vector space over the field $\KK$ 
and $G\subset\GL(V)$ be a finite subgroup. We
shall describe the Cox realization of the quotient space
$V/G:=\Spec\KK[V]^G$. Recall that a linear operator of finite
order $A\in\GL(V)$ is called {\it pseudoreflection}, if the subspace $V^A$
of $A$-fixed points is a hyperplane. By
Chevalley-Shephard-Todd's theorem, the algebra of invariants $\KK[V]^H$
of a finite subgroup $H\subset\GL(V)$ 
is free if and only if the subgroup $H$ is generated by pseudoreflections. Denote by $H$ the
subgroup generated by all pseudoreflections in $G$. It is
clear that $H$ is a normal subgroup of $G$. Let $\phi:G\to F:=G/H$
be the projection, $[F,F]$ be the commutant of the
group $F$, $N:=F/[F,F]$ and $\widetilde{H}:=\phi^{-1}([F,F])$. Put
$Z:=\Spec\KK[V]^{\widetilde{H}}$. The finite abelian group $N$
acts naturally on the variety $Z$.

\begin{theor} \label{tquot}
The quotient morphism $Z\xrightarrow{/N} V/G$ is the Cox
realization of the variety $V/G$.
\end{theor}

\begin{proof}
Denote the variety $\Spec\KK[V]^H$ by $W$. It is easy to prove that
the group $F$ acts on $W$ linearly. Let $\widetilde{F}\subseteq F$
be a subgroup generated by pseudoreflections.  Then the algebra of
invariants $\KK[W]^{\widetilde
F}\cong\KK[V]^{\phi^{-1}(\widetilde{F})}$ is free. This shows that
$F$ does not include any element acting on $W$ as a
pseudoreflection. Hence there exists an open $F$-invariant subset
$U\subseteq W$ such that $\codim_W(W\setminus U)\ge 2$ and the
$F$-action on $U$ is free. Let $\zeta:W\to W/[F,F]\cong Z$ be the
quotient morphism. Then $\codim_Z(Z\setminus\zeta(U))\ge 2$ and
the $N$-action on $\zeta(U)$ is free. According to
Proposition~\ref{prrel}, it is sufficient to prove that the $N$-variety
$Z$ is factorially graded. Let $D$ be an $N$-invariant Weil divisor on
$Z$ and $\zeta^{-1}(D)$ be its preimage in $W$. The divisor
$\zeta^{-1}(D)$ coincides with $\div(f)$, where $f$ is an
$F$-semiinvariant function on $W$. This implies that $f\in\KK[W]^{[F,F]}$, and
$D=\div(f)$ on $Z$.
\end{proof}

\begin{cor}
The variety $V/G$ is toric if and only if the group $G/H$ is
commutative.
\end{cor}

\begin{proof}
The Cox ring of a toric variety is a polynomial ring~\cite{Cox}.
Therefore the group $[F,F]$ is generated by pseudoreflections,
and hence it is trivial. The converse assertion follows from the
fact that a linear action of the finite abelian group $G/H$ on the
space $W$ is diagonazible.
\end{proof}

\begin{remark}
Assume that $R(X)$ is a polynomial ring. It is natural to ask
whether the variety $X$ is toric. In general, the answer is negative:
one may take a toric variety of dimension $\ge 2$ and remove a
finite set of points. Nevertheless it is true for a complete
variety with a free finitely generated divisor class group
\cite[Cor.~4.4]{BHC} The question whether this is true for affine
varieties reduces to the linearization problem for actions of
quasitori (see, for example, \cite{KR}), which is open.
\end{remark}

\begin{example} \label{exq8}
Let $V=\KK^2$ and
$$
G=Q_8=\left\{\pm E, \pm
\begin{pmatrix}
i & 0 \\
0 & -i
\end{pmatrix},
\pm
\begin{pmatrix}
0 & 1 \\
-1 & 0
\end{pmatrix},
\pm
\begin{pmatrix}
0 & i \\
i & 0
\end{pmatrix}\right\},
$$
where $i^2=-1$. Here $H=\{e\}$, $\widetilde{H}=\{\pm E\}$ and the
variety $Z$ is a two-dimensional quadratic cone. In particular,
the Cox ring of the variety $V/G$ is not factorial.
\end{example}

Theorem~\ref{tquot} can be partially generalized to the case of
infinite groups. With any affine algebraic group $G$ we associate
the group of characters $\mathfrak{X}(G)$ and the
intersection of kernels of all characters
$$G_1:=\cap_{\xi\in\mathfrak{X}(G)}\Ker(\xi).$$ It is clear that $G_1$
is a normal subgroup of $G$ and $N:=G/G_1$ is a quasitorus. If
$G$ is reductive, then $G_1$ is reductive too.

\begin{prop} \label{prrp}
Let $G$ be a reductive algebraic group and $Y$ be an irreducible
normal affine $G$-variety such that
\begin{enumerate}
\item $\KK[Y]^{\times}=\KK^{\times}$;
\item each $G$-stable Weil divisor on $Y$ is principle;
\item for the quotient morphism $\pi:Y\to Y/\!/G:=\Spec\KK[Y]^G$ there exists an open subset $U\subseteq Y/\!/G$ such
that $\codim_Y(Y\setminus\pi^{-1}(U))\ge 2$ and the $G$-action 
on $\pi^{-1}(U)$ is free.
\end{enumerate}
Set $Z=\Spec\KK[Y]^{G_1}$. Then the quotient morphism
$p:Z\xrightarrow{/\!/N} Y/\!/G$ is the Cox realization of
$Y/\!/G$.
\end{prop}

\begin{proof}
It is easy to see that $\codim_Z(Z\setminus p^{-1}(U))\ge 2$ and
the $N$-action on $p^{-1}(U)$ is free. Hence it is sufficient
to prove that each $N$-stable divisor $D$ on $Z$ is principle.
The preimage $\zeta^{-1}(D)$ of the divisor $D$ under the quotient
morphism $\zeta:Y\xrightarrow{/\!/G_1} Z$ is a $G$-stable divisor
on $Y$. Hence $\zeta^{-1}(D)=\div(f)$ for some $G$-semiinvariant
$f\in\KK[Y]$. Then $f\in\KK[Y]^{G_1}$, and $D=\div(f)$ on $Z$.
\end{proof}

\begin{example} (see~\cite[Sec.~3]{ArHa} and \cite[Th.~4.1]{Ar})
Let $\widehat{G}$ be a semisimple connected simply connected
affine algebraic group and $G\subseteq\widehat{G}$ be a reductive
subgroup. It is known that $\widehat{G}$ is a factorial variety,
$\KK[\widehat{G}]^{\times}=\KK^{\times}$ and the homogeneous
space $\widehat{G}/G$ is affine. Hence, the map
$\widehat{G}/G_1\xrightarrow{/N}\widehat{G}/G$ is the Cox
realization.
\end{example}

More generally, consider the Cox realization of the variety of
double cosets. Let $G$ and $G'$ be reductive subgroups
of a simply connected semisimple group $\widehat{G}$. Then $G\times G'$ 
acts on $\widehat{G}$ as
$(g,g')\hat{g}=g\hat{g}g'^{-1}$. Denote by
$G\backslash\widehat{G}/G'$ the categorical quotient for this
action. In order to prove that $G_1\backslash\widehat{G}/G'_1\to
G\backslash\widehat{G}/G'$ is the Cox realization, one needs to
check condition (iii) of Proposition~\ref{prrp}.

\begin{example}
Let $\widehat{G}=\SL(3)$, $G'=\SO(3)$ and
$$G=T^2=\{\diag(t_1,t_2,t_3): t_1t_2t_3=1\}.$$ Here
$Z=G_1\backslash\widehat{G}/G'_1=\SL(3)/\SO(3)$ is the variety of
symmetric matrices $(a_{ij})$, $i,j=1,2,3$,  with determinant equals 
1, $N=T^2$ acts on $Z$ as $a_{ij}\to t_it_ja_{ij}$,
and $X=G\backslash\widehat{G}/G'=Z/\!/T^2$ is a three-dimensional
hypersurface. As the subset $\pi^{-1}(U)$, one may take the
preimage in $\SL(3)$ of the set of symmetric matrices with at most one zero element 
$a_{ij}$, $i\ge j$. Thus the
spectrum of the Cox ring of the hypersurface $X$:
$$
2x_1^3+x_2x_3x_4-x_1^2-x_1^2x_2-x_1^2x_3-x_1^2x_4=0
$$
is the hypersurface $Z$:
$$
y_1y_2y_3+2y_4y_5y_6-y_1y_4^2-y_2y_5^2-y_3y_6^2=1.
$$
\end{example}

\begin{example}
Let $\widehat{G}=\SL(4)$, $G'=\Sp(4)$ and
$$G=T^3=\{\diag(t_1,t_2,t_3,t_4): t_1t_2t_3t_4=1\}.$$ Here
$Z=\SL(4)/\Sp(4)$ is the variety of skew-symmetric matrices
$(b_{ij})$, $i,j=1,2,3,4$, with Pfaffian
$\Pf=b_{12}b_{34}-b_{13}b_{24}+b_{14}b_{23}$ equals 1. 
In this case, condition (iii) of Proposition~\ref{prrp} is not fulfilled:
the closure of the $T^3$-orbit of any point on the divisor $\{b_{12}=0\}$
contains a point with $b_{34}=0$. It is easy to prove that
$X=T^3\backslash\SL(4)/\Sp(4)\cong\KK^2$. Therefore the spectrum
of the Cox ring can not be realized as a homogeneous space of the group $\SL(4)$.
\end{example}

\section{A universal property of the Cox realization}
\label{sec2} Let $Z$ be an irreducible affine variety with an action of
a quasitorus $Q$. We call the action {\it
uncontracting} if for any prime divisor $D\subset Z$ the
closure of its image $\overline{\pi(D)}$ under the quotient
morphism $\pi:Z\to Z/\!/Q$ has codimension $\leq1$ in $Z/\!/Q$.
Note that if $Q$ is a finite abelian group, then every $Q$-action
is uncontracting.

\begin{theor} \label{tuniv}
Let $Z$ be an irreducible normal affine variety with an uncontracting
action of a quasitorus $Q$ and $\pi:Z\to Z/\!/Q$ be the quotient morphism. 
Suppose that the $Q$-variety $Z$ is factorially graded and the
quotient space $X=Z/\!/Q$ has the Cox realization $q:\overline{X}\to X$. Then
there exists a surjective homomorphism $\mu:Q\to N$ and a dominant
morphism $\nu:Z\to\overline{X}$ such that $\nu(gz)=\mu(g)\nu(z)$
for all $g\in Q$, $z\in Z$, and the following diagram is
commutative:
$$
\xymatrix{ & Z \ar[dr]^{\pi} \ar[rr]^{\nu} & & \overline{X}
\ar[dl]_{q} &
\\
& & X. & & }
$$
Moreover, the homomorphism $\mu$ is unique,
and the morphism $\nu$ is determined up to composition with the
automorphism of $Z$ defined by an element of $Q$.
\end{theor}

This result can be proved using restriction of divisors to the
subset $X^{\text{\reg}}$ and lifting of the corresponding line
bundles to $Z$, see the proof of Proposition~\ref{prrel}. We leave
details to the reader.

Theorem~\ref{tuniv} can be reformulated in terms of corresponding
semigroups of association classes of homogeneous elements. This
leads to Theorem~\ref{dtheo}, which is of independent
interest. The proof of Theorem~\ref{dtheo} is given in
Section~\ref{secd}. On the other hand, it is shown below that
Theorem~\ref{dtheo} implies Theorem~\ref{tuniv}.

\begin{proof} \ {\it Existence.}\
Let $A=\KK[Z]$ and $A=\oplus_{u\in M}A_u$ be the grading
defined by the $Q$-action. By assumption, there is
an isomorphism $\delta:\KK[X]\to A_0$. We have to prove that there
exists a homogeneous embedding $\psi: R(X)\to A$
coinciding with $\delta$ on zero component and mapping
different homogeneous components of $R(X)$ to different components
of $A$. Let $S=\As(A)$ be the semigroup of association
classes of non-zero homogeneous elements of $A$,
$\Gamma=\As(\KK[X])$ and $D=\WDiv(X)_+$. The embedding $\tau:\Gamma\hookrightarrow D$ 
as the semigroup of effective divisors defines the
divisor theory of the semigroup $\Gamma$. The isomorphism
$\delta$ defines an embedding $\alpha:\Gamma\hookrightarrow S$ satisfying conditions $(*)$ and $(**)$ of
Theorem ~\ref{dtheo}. Indeed, the first condition follows from
existence of gradings and the second one is satisfied because the $Q$-action
is uncontracting. Thus there is an embedding
$\beta$ of the semigroup $D$, which is identified with $\As(R(X))$,
into the semigroup $S$ extending the embedding $\alpha$. We shall
"lift" this embedding to a homogeneous embedding of algebras.

Suppose the group $\Cl(X)$ has rank $r$. Let
$C_1,\dots,C_r,B_1,\dots,B_s$ be Weil divisors on $X$ such that
their classes generate $\Cl(X)$ with relations
$d_1[B_1]=\dots=d_s[B_s]=0$, $d_i\in\NN$. Fix non-zero elements of
the Cox ring $f_1,\dots,f_r,h_1,\dots,h_s$, $f_i\in R(X)_{[C_i]}$,
$h_j\in R(X)_{[B_j]}$ and non-zero homogeneous elements
$a_1,\dots,a_r,a_1',\dots,a_s'$ of $A$ such that
the classes of $f_i$ (resp. $h_j$) in $D$ map to the classes of $a_i$
(resp. $a_j'$) in $S$ under the embedding $\beta$. Note that
$h_j^{d_j}=F_j\in\KK[X]$. Therefore, we need the isomorphism
$\delta$ to map $F_j$ to $(a_j')^{d_j}$. This can be achieved by
changing the elements $a_j'$ with proportional ones.

The isomorphism $\delta$ extends to an isomorphism of the quotient
fields $\delta:\KK(X)\to QA_0$. Consider an element $u\in\Cl(X)$,
$$
u=u_1[C_1]+\dots+u_r[C_r]+v_1[B_1]+\dots+v_s[B_s], \ \ u_i\in\ZZ,
\ v_j\in\{0,\dots,d_j-1\},
$$
and the element $f_u:=f_1^{u_1}\dots f_r^{u_r}h_1^{v_1}\dots
h_s^{v_s}\in (QR(X))_u$. For any $f\in R(X)_u$, we have
$f=\frac{f}{f_u}f_u$ and $\frac{f}{f_u}\in (QR(X))_0=\KK(X)$. Put
$$
\psi(f)=\delta(\frac{f}{f_u})\psi(f_u)=\delta(\frac{f}{f_u})a_1^{u_1}\dots
a_r^{u_r}(a_1')^{v_1}\dots (a_s')^{v_s}.
$$
It is easy to check that the map $\psi$ is linear on any
homogeneous component $R(X)_u$ and corresponds to the embedding
$\beta$ of semigroups. In particular, the image of $\psi$ lies in
$A$. The grading defines a homomorphism of the semigroup $D$
(resp. $S$) to the group $\Cl(X)$ (resp. $M$), and the
kernel of this homomorphism is the semigroup $\Gamma$ (resp.
$\alpha(\Gamma)$). The semigroup $D$ projects onto the group
$\Cl(X)$ surjectively. Extending the above homomorphisms to
homomorphisms of groups generated by the semigroups, we obtain an
injective homomorphism $\Cl(X)\to M$. This proves the assertion of
the theorem concerning homogeneous components.

We have to check that for $f\in R(X)_u$ and $f'\in R(X)_{u'}$
the condition $\psi(ff')=\psi(f)\psi(f')$ holds. It is sufficient to
prove this for $f=f_u$, $f'=f_{u'}$, where
it follows from $\psi(F_j)=\delta(F_j)$,  $j=1,\dots,s$.
Thus the restriction of $\psi$ to any component
$R(X)_u$ is linear, $\psi$ is multiplicative on homogeneous
elements, and one can extend $\psi$ to a
homogeneous embedding of algebras $\psi:R(X)\to A$ via distributivity.

\smallskip
{\it Uniqueness.}\ Theorem~\ref{dtheo} implies that the embedding
of semigroups $\beta:D\hookrightarrow S$ is unique. Hence two
embeddings $\psi$ and $\psi'$ of algebras can differ on 
$a_1,\dots,a_r,a_1',\dots,a_s'$ only by scalar multiples; moreover, the
multiple for $a_j'$ has to be a root of
unity of degree $d_j$. Since the group $\Cl(X)$ is embedded into
$M$, this corresponds to the action of an element of $Q$ on $A$.
\end{proof}

The condition of graded factoriality of $Z$ in Theorem~\ref{tuniv} is
essential: one may take $Z$ to be the quotient of
$\overline{X}$ by a proper subgroup $Q\subset N$.
The following examples show that the condition of uncontracting is
also essential.

\begin{example}
Let $Z=\KK^3$, $Q=\KK^{\times}$ and
$t(z_1,z_2,z_3)=(tz_1,tz_2,t^{-2}z_3)$. Here $Z/\!/Q$ is
isomorphic to a two-dimensional quadratic cone $X$, and the torus
$Q$ does not admit a surjective homomorphism onto the Neron-Severi
quasitorus $N$ consisting of two elements. An analogous example of a
variety $X$  with a free divisor class group gives the cone
$X=\{(x_1,x_2,x_3,x_4) : x_1^2x_2-x_3x_4=0\}$ realized as the
quotient $Z/\!/Q$, where $Z=\KK^5$, $Q=(\KK^{\times})^2$ and
$$
(t_1,t_2)(z_1,z_2,z_3,z_4,z_5)=(t_1z_1,t_2z_2,t_1^{-1}t_2z_3,t_1t_2^{-1}z_4,t_1^{-1}t_2^{-1}z_5).
$$
\end{example}

\begin{example}
Let $Z=\KK^5$, $Q=(\KK^{\times})^2$ and
$$
(t_1,t_2)(z_1,z_2,z_3,z_4,z_5)=(t_1z_1,t_1z_2,t_2z_3,t_2z_4,t_1^{-1}t_2^{-1}z_5).
$$
The quotient space $Z/\!/Q$ is a three-dimensional quadratic cone
$X$. Here $\overline{X}=\KK^4$, $N=\KK^{\times}$ and the action
is given by $s(y_1,y_2,y_3,y_4)= (sy_1,sy_2,s^{-1}y_3,s^{-1}y_4)$.
In this case, there is an embedding of the Cox ring
$\KK[y_1,y_2,y_3,y_4]$ into $\KK[Z]$ extending the embedding of
semigroups, but it is not unique: one may take either
$$
y_1\to z_1, \ y_2\to z_2, \ y_3\to z_3z_5, \ y_4\to z_4z_5,
$$
or 
$$
y_1\to z_1z_5, \ y_2\to z_2z_5, \ y_3\to z_3, \ y_4\to z_4.
$$
\end{example}

%
\section{Lifting of automorphisms}
\label{sec3} 
Let $R=\oplus_{u\in M} R_u$ be an algebra graded by
a finitely generated abelian group $M$. Define a subgroup
$\widetilde{\Aut}(R)$ of the automorphism group $\Aut(R)$ as
$$
\widetilde{\Aut}(R)=\{ \phi\in\Aut(R) \mid \exists
\phi_0\in\Aut(M) : \phi(R_u)=R_{\phi_0(u)} \ \forall\, u\in M \}.
$$
We say that elements of $\widetilde{\Aut}(R)$
normalize the grading of the algebra $R$.

\begin{theor} \label{taut}
Let $X$ be an irreducible normal affine variety with a finitely
generated divisor class group $\Cl(X)$ and the Neron-Severi
quasitorus $N$. Assume that $\KK[X]^{\times}=\KK^{\times}$. Then
there is the following exact sequence:
$$
1 \to N \xrightarrow{\alpha} \widetilde{\Aut}(R(X))
\xrightarrow{\beta} \Aut(X) \to 1.
$$
\end{theor}
\begin{proof}

Each automorphism $\phi\in\widetilde{\Aut}(R(X))$ induces an
automorphism of the algebra $R(X)_0=\KK[X]$, and this defines the map
$\beta$. The quasitorus $N$ acts by homogeneous automorphisms of 
$R(X)$, which are identical on $R(X)_0$. This defines
$\alpha$ and shows that the composition $\beta\circ\alpha$ maps
elements of $N$ to the identical automorphism of $X$.

\begin{lemma} \label{lem1p}
Let $u\in\Cl(X)$ be a non-zero element. Then the $R(X)_0$-module
$R(X)_u$ is not cyclic.
\end{lemma}

\begin{proof}
Suppose that there is $f\in R(X)_u$ such that any $g\in R(X)_u$ has the form $g=fh$ 
for some $h\in R(X)_0$. Let $D$ be a
Weil divisor from the fixed set of representatives corresponding to the class $u$. 
Then the divisor
$(D+\div(g))-(D+\div(f))=\div(h)$ is effective. Put
$D+\div(f)=\sum_i a_iD_i$, $a_i\in\NN$. There exists a rational
function $F\in\KK(X)$ with pole of order one on $D_1$.
Multiplying it by a suitable regular function, we get
$\div(F)=-D_1+\sum_j b_jD_j$, $b_j\in\NN$. Then $fF\in R(X)_u$,
but $\div(fF)-\div(f)$ is not effective.
\end{proof}

Let us show that $\Ker\beta=\Imm\alpha$. Suppose
$\phi\in\Ker\beta$ and $f$ is a prime homogeneous element of 
$R(X)$. If $f\in R(X)_0$, then $\phi(f)=f$. Let
$f\in R(X)_u$, $u\ne 0$, and $g\in R(X)_u$ be an element
not divisible by $f$. The automorphism $\phi$ induces an
automorphism of the quotient field $QR(X)$ identical on $QR(X)_0$.
Since $\frac{f}{g}\in QR(X)_0$, one gets $f\phi(g)=g\phi(f)$. Since
$R(X)$ is factorially graded, $f$ divides $\phi(f)$. The element
$\phi(f)$ is also prime, therefore $\phi(f)=\lambda f$,
$\lambda\in\KK^{\times}$. Thus $\phi$ acts by scalar multiplication 
on prime homogeneous (and on all homogeneous)
elements. Moreover, $\phi\mid_{R(X)_u}$ is a scalar operator, so
$\phi$ defines a homomorphism from $\Cl(X)$ to $\KK^{\times}$,
that is, it belongs to $\alpha(N)$.

Now we have to prove that the map $\beta$ is surjective. Let
$\psi\in\Aut(X)$. Then $\psi$ induces automorphisms of the groups
$\WDiv(X)$ and $\Cl(X)$. Suppose $K$ is a subgroup of $\WDiv(X)$
from the definition of $R(X)$ (see Section~\ref{sec1}), and $\psi(K)$
is its image. Since the ring $R(X)$ does not depend (up to 
isomorphism) on the choices of $K$ and of a
coherent set $\Ff$, we can fix an isomorphism $\tau$ between
the Cox rings $R_{\psi(K),\psi^*(\Ff)}(X)$ and $R_{K,\Ff}(X)$ with the identical 
restriction to $R(X)_0=\KK[X]$. Then the
composition of $\psi^*: R_{K,\Ff}(X) \to
R_{\psi(K),\psi^*(\Ff)}(X)$ and $\tau$ is an element of the subgroup
$\widetilde{\Aut}(R(X))$, which induces 
the automorphism $\psi^*$ on $R(X)_0$.
\end{proof}

Let $X$ be an affine toric variety. 
The Cox ring $R(X)$ is a (graded) polynomial algebra with homogeneous generators \cite{Cox}.
Hence the description of the group $\Aut(X)$ may be reduced to the
description of the group of automorphisms of the polynomial algebra
normalizing the grading. In the study of automorphisms of the polynomial
algebra the concepts of tame and wild automorphisms play an important role.

\begin{de}
\begin{enumerate}
\item An automorphism of the algebra $\KK[y_1,\dots,y_m]$ is called {\it
elementary}, if it is either a linear map or a map of the
following form
$$
(y_1,\dots,y_m) \to (y_1,\dots,y_{i-1},y_i+f,y_{i+1},\dots,y_m),
$$
where $f\in\KK[y_1,\dots,y_{i-1},y_{i+1},\dots,y_m]$.
\item An automorphism is called {\it tame}, if it can be decomposed 
into a composition of elementary automorphisms.
\item An automorphism, which is not tame, is called {\it wild}.
\end{enumerate}
\end{de}

Define an {\it elementary} automorphism of the algebra $\KK[X]$
as the image under the map $\beta$ of an elementary automorphism of
the algebra $R(X)$ belonging to the subgroup
$\widetilde{\Aut}(R(X))$. Accordingly, we call an
automorphism {\it tame} if it can be represented as a composition
of elementary automorphisms, and {\it wild} otherwise. In the next
section we illustrate these concepts with the example of quadratic cone.

Recall that for any affine variety $X$ the group $\Aut(X)$ has a
structure of infinite-dimensional affine algebraic group in the sense
of~\cite{Sh}. This group is said to be {\it finite-dimensional} 
if $\Aut(X)$ has a structure of an affine algebraic group such that the
action $\Aut(X)\times X\to X$ is algebraic (an equivalent
definition of a finite-dimensional group of automorphisms is
given in~\cite{Ra}).

\begin{theor} \label{tinf}
Let $X$ be an affine toric variety of dimension $\ge 2$ with
$\KK[X]^{\times}=\KK^{\times}$. Then the group $\Aut(X)$ is not
finite-dimensional.
\end{theor}

\begin{proof}
Let $R(X)=\KK[y_1,\dots,y_m]$ be the Cox ring with the 
$\Cl(X)$-grading. The component $R(X)_0$ contains a monomial $h$,
which does not depend on $y_1$. Indeed, if it is not the case,
the Cox realization $q\colon \KK^m\rightarrow X$ contracts
the divisor $y_1=0$. But $y_1$ corresponds to a Weil divisor $D$
on $X$. It is easy to check that
$\pi(\div(y_1))=\mathrm{supp}(D)$. Lemma~\ref{lem1p} implies that
the homogeneous component containing $y_1$ contains a monomial
$f$, which does not depend on $y_1$. Consider homogeneous
automorphisms of $R(X)$ mapping $y_1$ to $y_1+fF(h)$ and
preserving the variables $y_2,\dots,y_m$, where $F(t)\in\KK[t]$. They
induce homomorphisms of the algebra $R(X)_0$. The images of any
monomial depending on $y_1$ are not contained in 
a finite-dimensional subspace. Thus the group $\Aut(X)$ is not finite-dimensional.
\end{proof}

%
\section{A wild automorphism of the quadratic cone}
\label{sec4} 

Consider the quadratic cone 
$$
X=\{(x_1,x_2,x_3,x_4) : x_1x_4-x_2x_3=0\}. 
$$
The variety $X$ can be realized as the cone
of degenerate $2\times 2$-matrices. Here $R(X)$ is the polynomial
algebra $\KK[y_1,y_2,y_3,y_4]$ with $\ZZ$-grading
$\deg(y_1)=\deg(y_2)=1$, $\deg(y_3)=\deg(y_4)=-1$, and
$x_1=y_1y_3$, $x_2=y_1y_4$, $x_3=y_2y_3$, $x_4=y_2y_4$. Consider
an automorphism $\tau$ of the cone $X$ defined by the formula
$$
\tau\colon
\begin{pmatrix}
x_1 & x_2 \\
x_3 & x_4
\end{pmatrix}
\mapsto
\begin{pmatrix}
x_1 & x_2+x_1(x_3-x_2) \\
x_3+x_1(x_3-x_2) & x_4+(x_3+x_2)(x_3-x_2)+x_1(x_3-x_2)^2
\end{pmatrix}.
$$
Its inverse is
$$
\tau^{-1}\colon
\begin{pmatrix}
x_1 & x_2 \\
x_3 & x_4
\end{pmatrix}
\mapsto
\begin{pmatrix}
x_1 & x_2-x_1(x_3-x_2) \\
x_3-x_1(x_3-x_2) & x_4-(x_3+x_2)(x_3-x_2)+x_1(x_3-x_2)^2
\end{pmatrix}.
$$
The lift of this automorphism to the Cox ring is defined as
$$
\zeta\colon(y_1,y_2,y_3,y_4)\mapsto(y_1,y_2+y_1(y_1y_4-y_2y_3),y_3,y_4+y_3(y_1y_4-y_2y_3)).
$$
This is known as {\it Anick's automorphism} of the polynomial algebra in
four variables.

\begin{theor} \label{twild}
The automorphism $\tau$ is wild.
\end{theor}

\begin{proof}
An elementary automorphism $\phi$ of $R(X)$ that preserves the
grading either multiplies the variables by non-zero scalars
or maps $(y_1,y_2,y_3,y_4)$ to one of the sets
$$
(y_1+y_2H_1(y_2y_3,y_2y_4),y_2,y_3,y_4), \
(y_1,y_2+y_1H_2(y_1y_3,y_1y_4),y_3,y_4),
$$
$$
(y_1,y_2,y_3+y_4H_3(y_1y_4,y_2y_4), y_4), \
(y_1,y_2,y_3,y_4+y_3H_4(y_1y_3,y_2y_3)).
$$
Therefore, the group of tame automorphisms is generated by these
automorphisms and the transpose automorphism $(y_1,y_2,y_3,y_4)\to
(y_3,y_4,y_1,y_2)$, which reverses the grading. Suppose
Anick's automorphism is equal to a composition of elementary
automorphisms normalizing the grading:
$$
\zeta=\phi_n\circ \dots \circ\phi_2\circ \phi_1.
$$
Since $\zeta$ preserves the grading and the result of conjugation
of an elementary automorphism preserving the grading with the transpose automorphism is an
elementary automorphism preserving the grading, one may assume
that each $\phi_i$ preserves the grading.

Let us replace all nonlinear automorphisms changing $y_3$ or $y_4$
in the above composition by their linear parts. Denote the new
composition by $\rho$. Since the linear parts of $\zeta $ and the identical
automorphism are the same, we obtain $\rho(y_3)=y_3$,
$\rho(y_4)=y_4$. Let $f=\zeta(y_1)-\rho(y_1)$ and
$g=\zeta(y_2)-\rho(y_2)$. Then
$$
\rho: (y_1,y_2,y_3,y_4) \to (y_1-f, y_2+y_1(y_1y_4-y_2y_3)-g, y_3,
y_4).
$$
Consider an ideal $I=(y_1,y_2)\lhd\KK[y_1,y_2,y_3,y_4]$.

\begin{lemma}
The polynomials $f$ and $g$ are in $I^3$.
\end{lemma}

\begin{proof}
Let $\zeta_j=\phi_j\circ \dots \circ\phi_1$ and $\rho_j$ be
the composition obtained from $\zeta_j$ by
replacing of automorphisms changing $y_3$ or $y_4$ with their linear
parts. Set $f_j=\zeta_j(y_1)-\rho_j(y_1)$,
$g_j=\zeta_j(y_2)-\rho_j(y_2)$, $h_j=\zeta_j(y_3)-\rho_j(y_3)$,
$s_j=\zeta_j(y_4)-\rho_j(y_4)$.

Let us prove by induction on $j$ that $f_j, g_j\in I^3$ and $h_j,
s_j\in I$. The case $j=1$ is easy.
At the step from $j=k$ to $j=k+1$, only one of the
polynomials $f_j,\ g_j,\ h_j,\ s_j$ is changed. We consider
the cases when $f_j$ and $h_j$ are changed, other cases are similar.

\smallskip
1) Assume that $f_j$ is changed. Then $\phi_{k+1}$ has the form
$$
(y_1,y_2,y_3,y_4)\mapsto(y_1+y_2H_1(y_2y_3,y_2y_4),y_2,y_3,y_4).
$$
We have
\begin{multline*}
f_{k+1}=\zeta_{k+1}(y_1)-\rho_{k+1}(y_1)= \zeta_k(y_1)+
\zeta_k(y_2)H_1(\zeta_k(y_2)\zeta_k(y_3),\zeta_k(y_2)\zeta_k(y_4))-\\-
\rho_k(y_1)-\rho_k(y_2)H_1(\rho_k(y_2)\rho_k(y_3),\rho_k(y_2)\rho_k(y_4))=
f_k+g_kH_1(\zeta_k(y_2)\zeta_k(y_3),\zeta_k(y_2)\zeta_k(y_4))+\\+
\rho_k(y_2)[H_1(\zeta_k(y_2)\zeta_k(y_3),\zeta_k(y_2)\zeta_k(y_4))-
H_1(\rho_k(y_2)\rho_k(y_3),\rho_k(y_2)\rho_k(y_4))].
\end{multline*}
By the inductive assumption, $f_k$ and $g_k$ belong to $I^3$. We have to
prove that the last summand belongs to $I^3$. It is sufficient to
do it when $H_1$ is a monomial. Let $H_1(u,v)=u^lv^r$.
The automorphisms $\zeta_k$ and $\rho_k$ preserve the grading, hence
$\zeta_k(y_2)\in I$ and $\rho_k(y_2)\in I$. Therefore, if $l+r\geq
2$, then
$$
\rho_k(y_2)[H_1(\zeta_k(y_2)\zeta_k(y_3),\zeta_k(y_2)\zeta_k(y_4))-
H_1(\rho_k(y_2)\rho_k(y_3),\rho_k(y_2)\rho_k(y_4))]\in I^3.
$$
If $l+r<2$, then either $l+r=0$, thus $H_1=\mathrm{const}$ and
$$
H_1(\zeta_k(y_2)\zeta_k(y_3),\zeta_k(y_2)\zeta_k(y_4))-
H_1(\rho_k(y_2)\rho_k(y_3),\rho_k(y_2)\rho_k(y_4))=0,
$$
or $l+r=1$, and without loss of generality one may assume that
$l=1,\ r=0$. One gets
\begin{multline*}
\rho_k(y_2)(H_1(\zeta_k(y_2)\zeta_k(y_3),\zeta_k(y_2)\zeta_k(y_4))-
H_1(\rho_k(y_2)\rho_k(y_3),\rho_k(y_2)\rho_k(y_4)))=\\
=\rho_k(y_2)(\zeta_k(y_2)\zeta_k(y_3)-\rho_k(y_2)\rho_k(y_3))
=\rho_k(y_2)((\rho_k(y_2)+g_k)(\rho_k(y_3)+h_k)-\rho_k(y_2)\rho_k(y_3))=\\
=\rho_k(y_2)(\rho_k(y_2)h_k+g_k\rho_k(y_3)+g_kh_k).
\end{multline*}
Since $\rho_k(y_2),h_k\in I$ and $g_k\in I^3$, 
the sum belongs to $I^3$. Thus, $f_{k+1}\in I^3$.

\smallskip
2) Suppose $h_j$ is changed. Then $\phi_{k+1}$ has the form
$$
(y_1,y_2,y_3,y_4)\mapsto(y_1,y_2,y_3+\mu
y_4+y_4F(y_1y_4,y_2y_4),y_4),
$$
where $F$ is a polynomial without constant term. Hence,
\begin{multline*}
h_{k+1}=\zeta_{k+1}(y_3)-\rho_{k+1}(y_3)=\\=
\zeta_k(y_3)+\mu\zeta_k(y_4)+
\zeta_k(y_4)F(\zeta_k(y_1)\zeta_k(y_4),\zeta_k(y_2)\zeta_k(y_4))-\rho_k(y_3)-
\mu\rho_k(y_4)=\\=h_k+\mu s_k+
\zeta_k(y_4)F(\zeta_k(y_1)\zeta_k(y_4),\zeta_k(y_2)\zeta_k(y_4)).
\end{multline*}
Since $h_k\in I$, $s_k\in I$, $\zeta_k(y_1)\in I$ and
$\zeta_k(y_2)\in I$, we obtain $h_{k+1}\in I$.
\end{proof}
Let us calculate the Jacobian matrix $J$ of the automorphism
$\rho$. Since a partial derivative of a polynomial from $I^3$
is in $I^2$, we obtain
$$
J=
\begin{pmatrix}
1 & 0 & 0 & 0 \\
2y_1y_4-y_2y_3 & 1-y_1y_3 & -y_1y_2 & y_1^2 \\
0 & 0 & 1 & 0 \\
0 & 0 & 0 & 1
\end{pmatrix}+ Q,
$$
where $Q$ is a matrix with elements in $I^2$.
Therefore,
$$
\det(J)=1-y_1y_3+c,\ c\in I^2.
$$
Since the determinant is not a constant, $\rho$ is not
an automorphism. This contradiction completes the proof of the theorem.
\end{proof}

It is interesting to note that if one adds the fifth $\zeta$-stable
variable of degree zero, the automorphism $\zeta$ becomes tame in
the class of automorphisms preserving the grading~\cite{Sm}. Let us
write down a decomposition of this automorphism into elementary ones:
\begin{multline*}
(y_1,y_2,y_3,y_4,y_5)\mapsto(y_1,y_2,y_3,y_4,y_5+(y_1y_4-y_2y_3))\mapsto\\(y_1,y_2+y_1y_5+y_1(y_1y_4-y_2y_3),y_3,y_4,y_5+(y_1y_4-y_2y_3))\mapsto\\
(y_1,y_2+y_1y_5+y_1(y_1y_4-y_2y_3),y_3,y_4+y_3y_5+y_3(y_1y_4-y_2y_3),y_5+(y_1y_4-y_2y_3))\mapsto\\
(y_1,y_2+y_1y_5+y_1(y_1y_4-y_2y_3),y_3,y_4+y_3y_5+y_3(y_1y_4-y_2y_3),y_5)\mapsto\\(y_1,y_2+y_1(y_1y_4-y_2y_3),y_3,y_4+y_3y_5+y_3(y_1y_4-y_2y_3),y_5)
\mapsto\\
(y_1,y_2+y_1(y_1y_4-y_2y_3),y_3,y_4+y_3(y_1y_4-y_2y_3),y_5).
\end{multline*}
\begin{remark}

The famous Nagata's automorphism
$$
(y_1,y_2,y_3) \to (y_1-2y_2(y_1y_3+y_2^2)-y_3(y_1y_3+y_2^2)^2,
y_2+y_3(y_1y_3+y_2^2), y_3)
$$
is a wild automorphism of the polynomial algebra in three
variables~\cite{SU}. This automorphism is homogeneous with respect to
a grading if and only if $\deg(y_1)=2\deg(y_2)$ and $\deg(y_3)=-\deg(y_2)$. 
For a $\ZZ$-grading, one gets the action 
$$
t\cdot(y_1,y_2,y_3)=(t^3y_1,ty_2,t^{-1}y_3)
$$
of a one-dimensional torus with a quotient space isomorphic to $\KK^2$.
Here the quotient morphism is not the Cox realization of the quotient space.
Nevertheless, in the case of a $\ZZ_n$-grading, the quotient morphism for the action  
$(\epsilon^3y_1,\epsilon y_2,\epsilon^{-1}y_3)$, $\epsilon^n=1$, is the Cox 
realization of $X_n:=\KK^3/\ZZ_n$, and Nagata's automorphism induces a wild 
automorphism of the variety $X_n$. In particular, Nagata's automorphism defines
a wild automorphism of the cone $X_2\subset\KK^6$ over the image
of the Veronese map $\PP^2\hookrightarrow\PP^5$.
\end{remark}

\medskip

%
\section{Appendix.  The divisor theory of a semigroup}
\label{secd} 
In this section a semigroup is a commutative
semigroup with unit $e$ such that all non-unit elements are not invertible.
We use multiplicative notation for the operation. A
semigroup $S$ is called {\it free}, if there is a subset $P\subset
S\setminus\{e\}$ such that each $s\in S\setminus\{e\}$ can be
expressed as $s=p_1^{k_1}\dots p_m^{k_m}$, $p_i\in P$, $k_i\in\NN$,
and this representation is unique up to the order of factors.
Elements of $P$ are prime (or indecomposable) elements of
$S$, so $P$ is uniquely determined by $S$. Sometimes it is
convenient to say "factorial semigroup" instead of "free semigroup".

Let $\Gamma$ be an arbitrary semigroup. The following definition
can be found in \cite[Ch.~III, \S~3]{BS}, see also \cite{Sk}.
\begin{de}
{\it The divisor theory} of a semigroup $\Gamma$ is an embedding
$\tau:\Gamma\hookrightarrow D$ into a free semigroup $D$ such that
\begin{enumerate}
\item
if $a,b\in\Gamma$ and $\tau(a)=\tau(b)c_1$ for some $c_1\in D$,
then $c_1=\tau(c)$ for some $c\in\Gamma$;
\item
if for $d_1,d_2\in D$ one has
$$
\{a\in\Gamma \mid \exists \, d\in D : \tau(a)=dd_1\}=\{b\in \Gamma
\mid \exists \, d\in D : \tau(b)=dd_2\},
$$
then $d_1=d_2$.
\end{enumerate}
\end{de}
In \cite[Ch.~III,\S~3,Th.~1]{BS}, it is proved that if a
semigroup $\Gamma$ has a divisor theory, then this theory is
unique up to isomorphism. So we denote $D$ as $D(\Gamma)$. Also we
do not distinguish between elements of $\Gamma$ and their
images in $D(\Gamma)$.

Suppose $\Gamma$ is a finitely generated semigroup. Let us explain
when $\Gamma$ admits the divisor theory, and recall a known
realization of this theory. Since there are only finitely many
prime elements of $D(\Gamma)$ in the decomposition of any
generating element of $\Gamma$, the semigroup $D(\Gamma)$ is
finitely generated and $\Gamma$ can be embedded into the lattice
$\ZZ D(\Gamma)$. The saturation condition emerging in the
following lemma is necessary for existence of the divisor theory of
$\Gamma$, see~\cite[Ch.~III,\S~3,Th.~3]{BS}.

\begin{lemma}
Let $\Gamma\subseteq D(\Gamma)$ be the divisor theory of a
semigroup $\Gamma$ and $L$ be the subgroup of the free abelian
group $\ZZ D(\Gamma)$ generated by $\Gamma$. If for an element
$l\in L$ there is $m\in\NN$ such that $l^m\in\Gamma$, then $l\in\Gamma$.
\end{lemma}

\begin{proof}
The condition $l^m\in\Gamma$ implies $l\in D(\Gamma)$. Since
$l=ab^{-1}$, $a,b\in\Gamma$, we have $a=bl$, and 
condition (ii) implies $l\in\Gamma$.
\end{proof}

It is known that for a finitely generated $\Gamma$ the saturation
condition is equivalent to the condition $\Gamma=\sigma\cap L$,
where $\sigma$ is the cone in the space
$L_{\QQ}:=L\otimes_{\ZZ}\QQ$ generated by (generators of)
$\Gamma$. Let $\sigma^{\vee}$ be the dual cone in the dual space
$L_{\QQ}^*$, that is
$$
\sigma^{\vee}=\{ f\in L_{\QQ}^* : f(x)\ge 0 \ \forall \ x\in\sigma
\},
$$
and let $f_1,\dots,f_r$ be primitive vectors on the edges of the cone
$\sigma^{\vee}$. Then the embedding
$$
\Gamma\hookrightarrow\ZZ^r_{\ge 0}, \ \ \ \ \ \ a \to
(f_1(a),\dots,f_r(a))
$$
is the divisor theory of $\Gamma$.
Conversely, let $P\subseteq\mathbb{Q}^r$ be a subspace such that
for the cone $\sigma:=P\cap\mathbb{Q}^r_{\geq 0}$ the subsets
$\sigma_i:=\sigma\cup\{x_i=0\}, i=1,\ldots,r$, are pairwise distinct 
and form the set of facets of the cones $\sigma$. Suppose also that
$L\subseteq P\cap\mathbb{Z}^r$ is a sublattice whose images under the coordinate
projections to the axis $Ox_i$ coincide with the set of
integer points on the $Ox_i$. Then the embedding
$\Gamma\hookrightarrow\mathbb{Z}^r_{\geq 0}$ is the divisor theory of
the semigroup $\Gamma:=L\cap\sigma$. Indeed, the set of coordinate
functions $x_i$ coincides with the set of primitive vectors on the edges of
$\sigma^{\vee}$.

\smallskip

Now we come to non-finitely generated semigroups. 
The theory of Weil divisors provides an example of
the divisor theory in algebraic geometry. Let $X$ be a normal
irreducible affine algebraic variety
and $\Gamma$ be the semigroup of association classes of the
algebra $\KK[X]$. The
semigroup $\Gamma$ can be identified with the semigroup of principle
effective divisors on $X$, and the semigroup $D(\Gamma)$ can be
realized as the semigroup of effective Weil divisors on $X$.
Examples of divisor theories appearing in Number Theory can be found
in \cite[Ch.~III]{BS}.

\smallskip

The aim of this appendix is to prove the following result
generalizing the uniqueness theorem for the divisor theory.

\begin{theor} \label{dtheo}
Let $\tau:\Gamma\to D(\Gamma)$ be the divisor theory of a
semigroup $\Gamma$ and $\alpha:\Gamma\to S$ be an embedding of
$\Gamma$ into a free semigroup $S$ satisfying the following
conditions:

 \ \ $(*)$ \ \ if for some $a,b\in\Gamma$ there exists $s\in S$ such
that $\alpha(a)=\alpha(b)s$, then $s=\alpha(c)$ for some
$c\in\Gamma$;

\ \ $(**)$\ \ if all elements of some subset $A\subseteq\Gamma$
are coprime in $D(\Gamma)$, that is, they are not divisible
by a non-unit element of $D(\Gamma)$, then all
elements of $\alpha(A)\subseteq S$ are coprime in $S$.

Then there exists a unique embedding $\beta:D(\Gamma)\hookrightarrow S$ such
that the following diagram is commutative:
$$
\xymatrix{ & \Gamma \ar[dr]^{\tau} \ar[rr]^{\alpha} & & S &
\\
& & D(\Gamma) \ar[ur]^{\beta}. & & }
$$
\end{theor}

\begin{proof} {\it Existence.}\
Let $P$ be the set of all prime elements of $S$ and $P_1\subseteq
P$ be the subset of elements that divide at least one of elements
$\alpha(a)$, $a\in\Gamma\setminus\{0\}$. Without loss of
generality we may assume that $P=P_1$. For any $d\in D:=D(\Gamma)$
and $s\in S$ we define
$$
L(d)=\{a\in\Gamma \mid \exists \, d'\in D : a=dd'\}; \ \
N(s)=\{b\in\Gamma \mid \exists \, s'\in S : \alpha(b)=ss'\}.
$$

\begin{lemma} \label{DL1}
\begin{enumerate}
\item Let $p\in P$. Then there is $q\in D$ such that
$L(q)\subseteq N(p)$.
\item Let $q\in S$ be a prime element. Then there is $p_1\in P$
with $N(p_1)\subseteq L(q)$.
\end{enumerate}
\end{lemma}

\begin{proof}
Since $N(p)\ne\emptyset$, there exists $a\in\Gamma$ such that
$\alpha(a)=ps$, $s\in S$. Let $a=q_1^{k_1}\dots q_l^{k_l}$ be the
prime decomposition in $D$. If
$L(q_i)\nsubseteq N(p)$, then for every $i=1,\dots,l$ there is
$a_i\in\Gamma$, which is divisible by $q_i$, and $\alpha(a_i)$ is
not divisible by $p$. Consider $b=a_1^{k_1}\dots a_l^{k_l}$. Then
$a$ divides $b$. But $\alpha(a)$ is divisible by $p$. It implies
that $p$ divides $\alpha(a_i)$, a contradiction.
Assertion (ii) can be proved in the same way.
\end{proof}

\begin{lemma} \label{DL2}
Let $q,q_1\in D$ be prime elements with $L(q)\subseteq L(q_1)$. Then $q=q_1$.
\end{lemma}

\begin{proof}
If the assertion of the lemma is false, then there is
$a\in\Gamma$, which is divisible by $q$ and is not divisible by
$qq_1$, hence it is in $L(q)\setminus L(q_1)$.
\end{proof}

For a fixed $q\in D$, there exist $p\in P$ and prime $q_1\in D$
such that $$L(q_1)\subseteq N(p)\subseteq L(q).$$ It implies
$q=q_1$ and $N(p)=L(q)$. Let $\{p_1,\dots,p_t\}$ be all prime
elements of $S$ satisfying $N(p_i)=L(q)$. Since $\cap_{j\in\NN}
N(p_i^j)=\emptyset$, for any $i$ there exists $r_i\in\NN$ such
that $$N(p_i)=N(p_i^2)=\dots=N(p_i^{r_i})\ne N(p_i^{r_i+1}).$$
Thus with any prime $q\in D$ one can associate
the set of prime elements $\{p_1,\dots,p_t\}$ and the set of exponents $\{r_1,\dots,r_t\}$.

Let us define $\beta(q)=p_1^{r_1}\dots p_t^{r_t}$ and extend this
map to a homomorphism $\beta:D\to S$. Since the sets
$\{p_1,\dots,p_t\}$ corresponding to different $q$ have empty
intersection, the map $\beta$ is injective. It remains to prove
that $\beta(a)=\alpha(a)$ for any $a\in\Gamma$.

Assume that a prime divisor $q\in D$ appears in the decomposition of
an element $a$ with multiplicity $k$, $\beta(q)=p_1^{r_1}\dots
p_t^{r_t}$, and any $p_i$ appears in the decomposition of
$\alpha(a)$ with multiplicity $n_i$. For every $p_i$, there exists
$q'\in D$ with $L(q')\subseteq N(p_i)$. If this inclusion is
strict, then there is an element $b\in N(p_i)$, which is not
divisible by $q'$. But then the elements of $L(q')$ and the element
$b$ are coprime in $D$. This contradicts to
condition $(**)$. Therefore, $L(q')=N(p_i)$. Hence each prime
divisor $p_i$ of the element $\alpha(a)$ corresponds to a (unique)
prime divisor $q$ of the element $a$, and it is sufficient to prove
that $n_i=kr_i$ for every $i=1,\dots,t$.

\begin{lemma} \label{DL3}
Suppose that $s\in S$ and $p_1,\dots,p_t$ do not divide $s$. Then there
exists $b\in\Gamma$ such that $\alpha(b)$ is divisible by $s$ and
is not divisible by $p_1,\dots,p_t$.
\end{lemma}

\begin{proof}
Let $s=s_1\dots s_n$ be the prime decomposition. If
$N(s_j)\subseteq N(p_i)=L(q)$, then, as we just have seen,
$N(s_j)=L(q)$, and $s_j$ coincides with one of the elements
$p_1,\dots,p_t$, a contradiction. Hence there are $b_j\in\Gamma$
such that $\alpha(b_j)$ is divisible by $s_j$ and not divisible
by $p_1,\dots,p_t$. Put $b=b_1\dots b_t$.
\end{proof}

We have $L(q)=N(p_1^{r_1}\dots p_t^{r_t})\supsetneqq
N(p_j^{r_j+1})$. Therefore for every $j$ there exists
$a_j\in\Gamma$ such that $\alpha(a_j)$ is divisible by
$p_1^{r_1}\dots p_t^{r_t}$, but not divisible by $p_j^{r_j+1}$.
Then
$$
\alpha(a_j)=p_1^{r_{1j}}\dots p_t^{r_{tj}}h_j, \ \ r_{ij}\ge r_i,
\ \ r_{jj}=r_j,
$$
and the element $h_j$ is coprime with $p_1,\dots,p_t$.

\begin{lemma} \label{DL4}
$r_{ij}=r_i$ for all $i,j=1,\dots,t$.
\end{lemma}

\begin{proof}
Assume that there is a pair $(u,v)$ with $r_{uv}>r_u$.
Subtract the vector $(r_{1v},\dots,r_{tv})$ from the vector
$(r_{1u},\dots,r_{tu})$. At the $u$-th position we obtain a
negative number. Let $k_1$ be the first position in the difference
with a negative coordinate. Let us add the
vector $(r_{1k_1},\dots,r_{tk_1})$ to this difference. We 
repeat this operation till a vector $(z_1,\dots,z_t)$
with nonnegative coordinates is obtained. We have
$$
(z_1,\dots,z_t)+(r_{1v},\dots,r_{tv})=(r_{1u},\dots,r_{tu})+\sum_{i=1}^m
(r_{1k_i},\dots,r_{tk_i}).
$$
Note that $m\ge 1$ and $z_{k_m}<r_{k_m}$, because $z_{k_m}$ is
obtained from a negative number by summation with
$r_{k_mk_m}=r_{k_m}$. Put $c=a_ua_{k_1}\dots a_{k_m}$. Recall
that $\alpha(a_v)=p_1^{r_{1v}}\dots p^{r_{tv}}h_v$, and
Lemma~\ref{DL3} implies existence of $b\in\Gamma$ such that
$\alpha(b)$ is divisible by $h_v$ and not divisible by
$p_1,\dots,p_t$. Then $bc=a_vb'$ for some $b'\in\Gamma$. Therefore
$\alpha(b')=p_1^{z_1}\dots p_t^{z_t}f$, where $f$ is not divisible
by $p_1,\dots,p_t$. From the conditions $z_{k_m}<r_{k_m}$ and
$N(p_{k_m})=N(p_{k_m}^{r_m})$ it follows that $z_{k_m}=0$. But
$N(p_i)=N(p_{k_m})$, so all $z_i$ equal zero. On the other hand,
$$
z_v=r_{vu}+\sum_{i=1}^m r_{vk_i}-r_{vv}\ge\sum_{i=1}^m r_{vk_i}>0,
$$
a contradiction.
\end{proof}

\begin{lemma} \label{DL5}
There is an element $c\in\Gamma$, which is divisible by $q$, but
not divisible by $q^2$, and $\alpha(c)$ is divisible by
$p_i^{r_i}$, but not divisible by $p_i^{r_i+1}$, $i=1,\dots,t$.
\end{lemma}

\begin{proof}
The previous lemma implies that the element $a_1$ satisfies the last two
properties and is divisible by $q$. Let us prove that any
element $c$ with these properties is not divisible by $q^2$. Suppose
$c=q^rh$, where the element $h\in D$ is coprime with $q$,
and $r>1$. Fix an element $b\in\Gamma$, which is divisible by $qh$
and not divisible by $q^2$. Then $b^r=cg$, $g\in\Gamma$ and $g$
is coprime with $q$. Hence $\alpha(b)^r=\alpha(c)\alpha(g)$
and $\alpha(g)$ is not divisible by $p_1,\dots,p_t$. Therefore
$\alpha(b)$ is not divisible by $p_i^{r_i}$, since $r>1$. On the
other hand, $b$ is divisible by $q$, and hence $p_1^{r_1}\dots
p_t^{r_t}$ divides $\alpha(b)$, a contradiction.
\end{proof}

Let us return to the proof of the equality $n_i=kr_i$. One can
find $m$ such that $mr_i\le n_i$ for all $i=1,\dots,t$ and
$(m+1)r_j>n_j$ for at least one $j$. Then there exist
$f_1,f_2\in\Gamma$ such that $f_1a=f_2c^m$ and $\alpha(f_1)$ is
coprime with $p_1,\dots,p_t$. Multiplicity of $p_j$ in
the decomposition of $\alpha(f_2)$ is less then $r_j$, hence $f_2$
is coprime with $p_1,\dots,p_t$ and $f_2$ is coprime
with $q$. Thus $n_i=mr_i$. On the other hand, comparing
multiplicities of $q$ in the equation $f_1a=f_2c^m$, we
obtain $k=m$. This completes the proof of the equality $n_i=kr_i$.

\smallskip
{\it Uniqueness.}\ Let $\gamma:
D(\Gamma)\hookrightarrow S$ be another embedding satisfying the conditions of Theorem~\ref{dtheo}. 
Let $a\in\Gamma$, $a=q_1^{k_1}\dots q_l^{k_l}$. Then
$$
\alpha(a)=\beta(q_1)^{k_1}\dots\beta(q_l)^{k_l}=\gamma(q_1)^{k_1}\dots\gamma(q_l)^{k_l}.
\ \ \ \ \ \ \ \ \ \ \ \ \ \ \ \ \ \ \ \ \ \ \ \ \ (1)
$$
If $p$ is a prime factor in the decomposition of
$\gamma(q_i)$, then $L(q_i)\subseteq N(p)$, hence
$L(q_i)=N(p)$. Therefore $\gamma(q_i)=p_1^{m_1}\dots p_t^{m_t}$,
where $p_1,\dots,p_t$ are prime elements corresponding to $q_i$.
If some $m_j$ is greater than $r_j$, then one gets a contradiction
with $N(p_j^{r_j+1})\subsetneq L(q_i)$. Hence $m_j\le r_j$ for all
$j$, and (1) implies $m_j=r_j$. This proves that $\beta=\gamma$.
\end{proof}

Let us show that conditions $(*)$ and $(**)$ of Theorem~\ref{dtheo} are essential.

\begin{example}
Let $\Gamma$ be the subsemigroup of the multiplicative semigroup $\NN$
of positive integers generated by $10$, $14$, $15$ and $21$.
The divisor theory of $\Gamma$ is the embedding of $\Gamma$ into
the subsemigroup $D$ generated by  $2$, $3$, $5$ and $7$. On the
other hand, the embedding $\alpha:\Gamma\hookrightarrow \NN$ defined on
generators by $\alpha(10)=10$, $\alpha(14)=2$, $\alpha(15)=15$ and
$\alpha(21)=3$ can not be extended to an embedding of $D$ into $\NN$.
Here condition $(*)$ is not satisfied.
\end{example}

\begin{example}
Let $\Gamma$ be the subsemigroup of the multiplicative semigroup $\NN$
of positive integers generated by $4$, $6$ and $9$. The divisor theory
of $\Gamma$ is the embedding of $\Gamma$ into the
subsemigroup $D$ generated by $2$ and $3$. On the other hand, the
embedding $\alpha:\Gamma\hookrightarrow \NN$ defined on generators by
$\alpha(4)=20$, $\alpha(6)=30$ and $\alpha(9)=45$ can not be extended
to an embedding of $D$ into $\NN$. Here condition $(**)$ is not satisfied.
\end{example}


%
\end{document}